\documentclass[12pt]{article}

\usepackage{graphicx}

\usepackage{bm,hyperref}

\begin{document}
\title{Connections Between Real Polynomial Solutions Of Hypergeometric-type  
Differential Equations With Rodrigues Formula}
\author{H. J. Weber\\Department of Physics\\University of Virginia\\
Charlottesville, VA 22904, USA}
\date{\today}
\maketitle
\begin{abstract}
Starting from the Rodrigues representation of polynomial solutions of the 
general hypergeometric-type differential equation complementary polynomials 
are constructed using a natural method. Among the key results is a generating 
function in closed form leading to short and transparent derivations of 
recursion relations and addition theorem. The complementary polynomials 
satisfy a hypergeometric-type differential equation themselves, have a 
three-term recursion among others and obey Rodrigues formulas. Applications 
to the classical polynomials are given.
\end{abstract}

\vspace{3ex}
\leftline{MSC: 33C45, 34B24, 35Q40, 42C05}
\leftline{PACS codes: 02.30.Gp; 02.30.Hq; 02.30.Jr, 03.65.Ge}
\leftline{Keywords: Polynomials with Rodrigues formula; solutions of} 
\leftline{~~~~~~~~~~~~~~~hypergeometric-type differential equation;}
\leftline{~~~~~~~~~~~~~~~generating function in closed form; recursion} 
\leftline{~~~~~~~~~~~~~~~relations; addition theorem}
Running title: Connections between polynomials with Rodrigues

\section{Introduction}

Real polynomial solutions $P_l(x)$ of the hypergeometric-type differential 
equation (ODE)
\begin{eqnarray}
\sigma(x)\frac{d^2 P_l}{dx^2}+\tau(x)\frac{d P_l}{dx}+\Lambda_l P_l(x)=0,\ 
\Lambda_l=-l\tau'-\frac{l}{2}(l-1)\sigma''
\label{ode}
\end{eqnarray}
with $l=0, 1,\ldots$ and real, first and second-order coefficient polynomials 
\begin{eqnarray}
\sigma(x)=ex^2+2fx+g,\ \tau=a_l+b_l x
\label{sig}
\end{eqnarray}
are analyzed in ref.~\cite{ni}, \cite{is}. The (unnormalized) polynomials are 
generated from the Rodrigues formula
\begin{eqnarray}
P_l(x)=\frac{1}{w(x)}\frac{d^l}{dx^l}(\sigma^l(x)w(x)),~l=0, 1, \ldots,
\label{rod}
\end{eqnarray} 
where $w(x)$ is the possibly $l$ dependent weight function on the fundamental 
interval $(a,b)$ that satisfies Pearson's ODE
\begin{eqnarray}
\sigma(x)w'(x)=[\tau(x)-\sigma'(x)] w(x) 
\label{wode}
\end{eqnarray} 
to assure the self-adjointness of the differential operator of the 
hypergeometric ODE. Polynomial solutions of ODEs with $l$ dependent 
coefficients are studied in ref.~\cite{les} along with their orthogonality 
properties and zero distributions, which we therefore do not address here.    

Here our first goal is to construct complementary polynomials for them by 
reworking their Rodrigues representation, Eq.~(\ref{rod}), in a simple and 
natural way. The generating function of these complementary polynomials is 
obtained in closed form allowing for short and transparent derivations of  
general properties shared by the complementary polynomials. 
 
The paper is organized as follows. In the next section we introduce and 
construct the complementary polynomials. In Section~3 we establish 
their generating function, the key result from which recursion relations and an
addition theorem are derived in Section~4. The Sturm-Liouville ODE is derived 
in Section~5. Classical polynomial examples are given in Section~6.  

\section{Complementary Polynomials} 

{\bf Definition.} We now introduce the complementary polynomials 
${\cal P}_\nu(x;l)$ defining them in terms of the generalized Rodrigues 
representation 
\begin{eqnarray}
P_l(x)=\frac{1}{w(x)}\frac{d^{l-\nu}}{dx^{l-\nu}}\left(\sigma(x)^{l-\nu}w(x)
{\cal P}_\nu(x;l)\right), 
\label{grod}
\end{eqnarray}
where $\nu=0,1, \ldots, l;\ l=0, 1, \ldots.$    

{\bf Theorem 1.} ${\cal P}_\nu(x;l)$ is a polynomial of degree $\nu$ that 
satisfies the recursive differential equation:    
\begin{eqnarray}
{\cal P}_{\nu+1}(x;l)=\sigma(x)\frac{d{\cal P}_\nu(x;l)}{dx}
+[\tau(x)+(l-\nu-1)\sigma'(x)]{\cal P}_\nu(x;l).
\label{rode}
\end{eqnarray}
By the Rodrigues formula~(\ref{rod}), ${\cal P}_{0}(x;l)\equiv 1.$ 

{\bf Proof.} Equations~(\ref{grod}), and ~(\ref{rode}) follow by induction. 
The first step, $\nu=1$, is derived by carrying out explicitly the innermost 
differentiation in Eq.~(\ref{rod}), which is a natural way of working with the 
Rodrigues formula~(\ref{rod}) that yields 
\begin{eqnarray}
P_l(x)=\frac{1}{w(x)}\frac{d^{l-1}}{dx^{l-1}}\left(l\sigma^{l-1}(x)w(x)
\sigma'(x)+\sigma^l(x)w'(x)\right)
\label{s1}
\end{eqnarray}
showing, upon substituting Pearson's ODE~(\ref{wode}), that 
\begin{eqnarray}
{\cal P}_1(x;l)=(l-1)\sigma'(x)+\tau(x).
\label{ex1}
\end{eqnarray}
Assuming the validity of the Rodrigues formula~(\ref{grod}) for $\nu$ we carry 
out another differentiation in Eq.~(\ref{grod}) obtaining 
\begin{eqnarray}\nonumber
P_l(x)&=&\frac{1}{w(x)}\frac{d^{l-\nu-1}}{dx^{l-\nu-1}}\bigg\{(l-\nu)
\sigma(x)^{l-\nu-1}\sigma'(x)w(x){\cal P}_\nu(x;l)\\
\nonumber&+&\sigma^{l-\nu}(x)w'(x){\cal P}_\nu(x;l)+\sigma(x)^{l-\nu}w(x)
{\cal P'}_\nu(x;l)\bigg\}\\\nonumber&=&\frac{1}{w(x)}\frac{d^{l-\nu-1}}
{dx^{l-\nu-1}}(\sigma(x)^{l-\nu-1}w(x)[(l-\nu)\sigma'(x){\cal P}_{\nu}\\
\nonumber &+&(\tau-\sigma'(x)){\cal P}_\nu(x;l)+\sigma{\cal P'}_\nu(x;l)])\\
&=&\frac{1}{w(x)}\frac{d^{l-\nu-1}}{dx^{l-\nu-1}}\left(
\sigma(x)^{l-\nu-1}w(x){\cal P}_{\nu+1}(x;l)\right).
\label{s2}
\end{eqnarray}
Comparing the rhs of Eq.~(\ref{s2}) proves Eq.~(\ref{grod}) by induction along 
with the recursive ODE~(\ref{rode}) which allows constructing systematically 
the complementary polynomials starting from ${\cal P}_0(x;l)\equiv 1.$ For 
example, $\nu=0$ of the recursive ODE~(\ref{rode}) confirms Eq.~(\ref{ex1}).  

In terms of a generalized Rodrigues representation we have\\{\bf Theorem 2.} 
The polynomials ${\cal P}_{\nu}(x;l)$ satisfy the Rodrigues formulas
\begin{eqnarray}
{\cal P}_{\nu}(x;l)=w^{-1}(x) \sigma^{\nu-l}(x)\frac{d^{\nu}}{dx^{\nu}}
[w(x)\sigma^{l}(x)];
\label{rods}
\end{eqnarray}
\begin{eqnarray}
{\cal P}_{\nu}(x;l)=w^{-1}(x)\sigma^{\nu-l}(x)\frac{d^{\nu-\mu}}
{dx^{\nu-\mu}}\left(\sigma^{l-\mu}(x)w(x) {\cal P}_{\mu}(x;l)\right).
\label{b5}
\end{eqnarray}

{\bf Proof.} We prove the Rodrigues formulas for the ${\cal P}_{\nu}(x;l)$ 
polynomials by integrating first the homogeneous ODE~(\ref{rode}) while 
dropping the inhomogeneous term ${\cal P}_{\nu+1}(x;l)$. This yields
\begin{eqnarray}
\ln {\cal P}_{\nu}(x;l)-\ln c_{\nu}=(-l+\nu+1)\ln\sigma(x)
-a_l\int\frac{dx}{\sigma(x)}-b_l\int\frac{xdx}{\sigma(x)},
\label{b1}
\end{eqnarray}
where $c_{\nu}$ is an integration constant and $\int\frac{dx}{\sigma(x)}, 
\int\frac{xdx}{\sigma(x)}$ are indefinite integrals. Exponentiating 
Eq.~(\ref{b1}) we obtain
\begin{eqnarray}
{\cal P}_{\nu}(x;l)=c_{\nu}\sigma(x)^{-l+\nu+1} e^{-a_l\int\frac{dx}
{\sigma(x)}-b_l\int\frac{xdx}{\sigma(x)}}.
\label{b2}
\end{eqnarray}
Note that, if the zeros of $\sigma(x)$ are real, they lie outside the 
fundamental interval $(a,b)$ of $w(x)$ and the hypergeometric Eq.~(\ref{ode}) 
by definition, while $x$ lies within it. So, these zeros pose no problem for 
the indefinite integrals.   

Now we allow for the $x$ dependence of $c_{\nu}$ and vary it to include the
inhomogeneous term ${\cal P}_{\nu+1}(x;l)$. Differentiating Eq.~(\ref{b2})
and substituting the recursive ODE~(\ref{rode}) yields
\begin{eqnarray}
{\cal P}_{\nu+1}(x;l)=c'_{\nu}(x)\sigma^{-l+\nu+2}(x)e^{-a_l\int\frac{dx}
{\sigma(x)}-b_l\int\frac{xdx}{\sigma(x)}},
\label{b3}
\end{eqnarray}
or 
\begin{eqnarray}\nonumber
{\cal P}_{\nu+1}(x;l)\sigma^{l-2-\nu}(x)e^{a_l\int\frac{dx}{\sigma(x)}+\int
\frac{xdx}{\sigma(x)}}&=&\frac{d}{dx}[\sigma^{l-\nu-1}(x)e^{a_l\int\frac{dx}
{\sigma(x)}+b_l\int\frac{xdx}{\sigma(x)}}{\cal P}_{\nu}(x;l)]\\
&=&c'_{\nu}(x).
\label{b4}
\end{eqnarray}
Noting that the expression in brackets on the rhs of Eq.~(\ref{b4}) differs 
from the coefficient of ${\cal P}_{\nu+1}(x;l)$ on the lhs only by one unit in 
the exponent of $\sigma(x)$ suggests iterating the differentiation and then 
replacing $\nu+1\to \nu.$ This leads to the formula
\begin{eqnarray}
{\cal P}_{\nu}(x;l)=\sigma^{-l+1+\nu}e^{-a_l\int\frac{dx}{\sigma}
-b_l\int\frac{xdx}{\sigma}}\frac{d^\nu}{dx^\nu}[\sigma^{l-1}
e^{a_l\int^x\frac{dx}{\sigma}+b_l\int\frac{xdx}{\sigma}}].
\label{explic}
\end{eqnarray}
Integrating Pearson's ODE~(\ref{wode}),
\begin{eqnarray}
\ln w(x)=\int\left(\frac{\tau}{\sigma}-\frac{\sigma'}{\sigma}\right)dx
=-\ln\sigma(x)+a_l\int\frac{dx}{\sigma(x)}+b_l\int\frac{xdx}{\sigma(x)}
\end{eqnarray}
and exponentiating this, gives
\begin{eqnarray}
w(x)=\sigma^{-1}e^{a_l\int\frac{dx}{\sigma(x)}+b_l\int\frac{xdx}{\sigma(x)}}.
\end{eqnarray}
Substituting this result into Eq.~(\ref{explic}) allows casting it in the 
general form of Eq.~(\ref{rods}). 

When we carry out the innermost differentiation in Eq.~(\ref{rods}) we
obtain the first step ($\mu=1$) of the inductive proof of the 
{\it generalized} Rodrigues representation of Eq.~(\ref{b5}). 
Equation~(\ref{b5}) yields trivially ${\cal P}_{\nu}(x;l)$ for
$\mu=\nu$, while for $\mu=\nu-1$ it reproduces Eq.~(\ref{rode}) and the case
$\mu=1$ is Eq.~(\ref{rods}). The inductive step from $\mu$ to $\mu+1$ is
similar to that leading to Eqs.~(\ref{grod}) and (\ref{rode}). 

\section{Generating Function}

{\bf Definition.} The generating function for the polynomials 
${\cal P}_\nu(x;l)$ is 
\begin{eqnarray}
{\cal P}(y,x;l)=\sum_{\nu=0}^{\infty}\frac{y^{\nu}}{\nu!}{\cal P}_{\nu}(x;l). 
\label{genfdef}
\end{eqnarray}
The series converges for $|y|<\epsilon$ for some $\epsilon>0$ and can be 
summed in closed form if the generating function is regular at the point $x$.

{\bf Theorem 3.} The generating function for the polynomials 
${\cal P}_\nu(x;l)$ is given in closed form by
\begin{eqnarray}
{\cal P}(y,x;l)=\frac{w(x+y\sigma(x))}{w(x)}\left(\frac{\sigma(x+y\sigma(x))}
{\sigma(x)}\right)^l;
\label{3_8}
\end{eqnarray}
\begin{eqnarray}
\frac{\partial^\mu}{\partial y^\mu}{\cal P}(y,x;l)=\frac{w(x+y\sigma(x))}{w(x)}
\left(\frac{\sigma(x+y\sigma(x))}{\sigma(x)}\right)^{l-\mu}
{\cal P}_\mu(x+y\sigma(x);l).
\label{3_6}
\end{eqnarray}
{\bf Proof.} Equation~(\ref{3_8}) follows by substituting the Rodrigues 
representation, Eq.~(\ref{rods}) in Eq.~(\ref{genfdef}) which yields, with 
$z\equiv x+y\sigma(x)$,
\begin{eqnarray}\nonumber
{\cal P}(y,x;l)&=&\sum_{\nu=0}^{\infty}\frac{y^{\nu}}{\nu!}{\cal P}_{\nu}(x;l)
=\left(w(x)\sigma^{l}(x)\right)^{-1}\sum_{\nu=0}^{\infty}
\frac{(y\sigma(x))^{\nu}}{\nu!}\frac{d^{\nu}}{dx^{\nu}}\left(\sigma^{l}(x)w(x)
\right)\\&=&\left(w(x)\sigma^{l}(x)\right)^{-1}\sum_{\nu=0}^{\infty}
\frac{(z-x)^{\nu}}{\nu!}\frac{d^{\nu}}{dz^{\nu}}
\left(\sigma^{l}(z)w(z)\right)|_{z=x},
\label{3_7}
\end{eqnarray}
converging for $|y\sigma(x)|<\epsilon$ for a suitable $\epsilon>0$ if 
$x\in (a,b)$ is a regular point of the generating function, i.e. $w$ is 
regular at $x$ and $x+y\sigma(x).$ The series can be summed exactly 
because the expression inside the derivatives is independent of the summation 
index $\nu$ and we deal with the Taylor expansion of the function 
$\sigma^{l}(z)w(z)$ at the point $x$ with increment $y \sigma(x)$. 

Differentiating Eq.~(\ref{genfdef}) and substituting the generalized Rodrigues 
formula~(\ref{b5}) in this yields Eq.~(\ref{3_6}) similarly.   

In preparation for recursion relations we translate the case $\mu=1$ of 
Eq.~(\ref{3_6}) into partial differential equations (PDEs).\\   
{\bf Theorem 4.} The generating function satisfies the partial differential
equations (PDEs)
\begin{eqnarray}\nonumber
(1+y\sigma'(x)+\frac{1}{2}y^2\sigma''\sigma(x))\frac{\partial 
{\cal P}(y,x;l)}{\partial y}&=&[{\cal P}_1(x;l)+y\sigma(x){\cal P}_1'(x;l)]\\
&\cdot&{\cal P}(y,x;l),
\label{pdey1}
\end{eqnarray}
\begin{eqnarray}
\frac{\partial {\cal P}(y,x;l)}{\partial y}=
[(l-1)\sigma'(x+y\sigma(x))+\tau(x+y\sigma(x))]{\cal P}(y,x;l-1),
\label{pdey2}
\end{eqnarray}
\begin{eqnarray}\nonumber
&&\left(1+y\sigma'(x)+\frac{1}{2}y^2\sigma''\sigma(x)\right)
\frac{\partial {\cal P}(y,x;l)}{\partial x}
\\&=&{\cal P}(y,x;l)y\bigg\{(1+y\sigma'(x)){\cal P}_1'(x;l)-\frac{1}
{2}y\sigma''{\cal P}_1(x;l)\bigg\},
\label{pdex1}
\end{eqnarray}
\begin{eqnarray}\nonumber
\sigma(x)\frac{\partial {\cal P}(y,x;l)}{\partial x}&=&(1+y\sigma'(x)))
[\tau(x)+(l-1)\sigma'(x)+y\sigma(x)(\tau'+(l-1)\sigma'')]\\&\cdot&
{\cal P}(y,x;l-1)-[\tau(x)+(l-1)\sigma'(x)]{\cal P}(y,x;l).
\label{pdex2}
\end{eqnarray}
{\bf Proof.} From Eq.~(\ref{3_6}) for $\mu=1$ in conjunction with 
Eq.~(\ref{3_8}) we obtain 
\begin{eqnarray}\nonumber
\sigma(x+y\sigma(x))\frac{\partial {\cal P}(y,x;l)}{\partial y}&=&
\sigma(x)[\tau(x+y\sigma(x))\\ &+&(l-1)\sigma'(x+y\sigma(x))]{\cal P}(y,x;l).
\label{pdexs}
\end{eqnarray}
Substituting in Eq.~(\ref{pdexs}) the Taylor series-type expansions  
\begin{eqnarray}\nonumber
\sigma(x+y\sigma(x))&=&\sigma(x)(1+y\sigma'(x)+\frac{1}{2}y^2\sigma''\sigma(x))
,\\\nonumber \sigma'(x+y\sigma(x))&=&\sigma'(x)+y\sigma''\sigma(x),\\
\tau(x+y\sigma(x))&=&\tau(x)+y\tau'(x)\sigma(x)
\label{sis}
\end{eqnarray}
following from Eq.~(\ref{sig}), we verify Eq.~(\ref{pdey1}). Using the 
exponent $l-1$ instead of $l$ of the generating function we can similarly 
derive Eq.~(\ref{pdey2}).   

By differentiation of the generating function, Eq.~(\ref{3_7}), with respect 
to the variable $x$ we find Eq.~(\ref{pdex2}). Using the exponent $l$ instead 
of $l-1$ of the generating function in conjunction with Eq.~(\ref{pdex2}) 
leads to Eq.~(\ref{pdex1}).   

\section{Recursion and Other Relations}

Our next goal is to rewrite various PDEs for the generating function in terms
of recursions for the complementary polynomials.\\{\bf Theorem 5.} The 
polynomials ${\cal P}_\nu(x;l)$ satisfy the recursion relations
\begin{eqnarray}\nonumber
{\cal P}_{\nu+1}(x;l)&=&[\tau(x)+(l-1-\nu)\sigma'(x)]{\cal P}_{\nu}(x;l)\\
&+&\nu\sigma(x)[\tau'+(l-1-\frac{1}{2}(\nu-1))\sigma'']{\cal P}_{\nu-1}(x;l);
\label{recu1}\\\nonumber
{\cal P}_{\nu+1}(x;l)&=&[\tau(x)+(l-1)\sigma'(x)]{\cal P}_{\nu}(x;l-1)\\
\nonumber &+&\nu\sigma(x)[\tau'+(l-1)\sigma'']{\cal P}_{\nu-1}(x;l-1)\\
&=&{\cal P}_{1}(x;l){\cal P}_{\nu}(x;l-1)+{\cal P}_{1}'(x;l-1)
{\cal P}_{\nu-1}(x;l-1);
\label{recu2}\\\nonumber
&\nu&(\nu-1)\frac{1}{2}\sigma''\sigma(x)\frac{d{\cal P}_{\nu-2}(x;l)}{dx}+
\nu\sigma'(x)\frac{d{\cal P}_{\nu-1}(x;l)}{dx}+\frac{d{\cal P}_{\nu}(x;l)}{dx}
\\\nonumber&=& \nu[\tau'+(l-1)\sigma'']{\cal P}_{\nu-1}(x;l)+\nu(\nu-1)
{\cal P}_{\nu-2}(x;l)\\\nonumber&\cdot&\bigg\{\sigma'(x)[\tau'+(l-1)\sigma'']
-\frac{1}{2}\sigma''[\tau(x)+(l-1)\sigma'(x)]\bigg\}\\\nonumber &=&
\nu{\cal P}_{1}'(x;l){\cal P}_{\nu-1}(x;l)+\nu(\nu-1){\cal P}_{\nu-2}(x;l)\\
&\cdot&\bigg\{\sigma'(x){\cal P}_{1}'(x;l)-\frac{1}{2}\sigma''{\cal P}_{1}(x;l)
\bigg\}. 
\label{recu3}
\end{eqnarray}
{\bf Proof.} Substituting Eq.~(\ref{genfdef}) defining the generating function 
in Eq.~(\ref{pdey1}) we rewrite the PDE as Eq.~(\ref{recu1}). The 
recursion~(\ref{recu2}) is derived similarly from Eq.~(\ref{pdey2}). The same 
way Eq.~(\ref{pdex1}) translates into the differential recursion 
relation~(\ref{recu3}). 

{\bf Corollary.} 
Comparing the recursion~(\ref{recu1}) with the recursive ODE~(\ref{rode}) we 
establish the basic recursive ODE
\begin{eqnarray}
\frac{d}{dx}{\cal P}_{\nu}(x;l)=\nu[\tau'+(l-1-\frac{1}{2}(\nu-1)\sigma'')]
{\cal P}_{\nu-1}(x;l)
\label{brecu}
\end{eqnarray}
with a coefficient that is independent of the variable $x$. 

{\bf Parameter Addition Theorem.}
\begin{eqnarray}
{\cal P}(y,x;l_1+l_2){\cal P}(y,x;0)={\cal P}(y,x;l_1){\cal P}(y,x;l_2). 
\label{addid}
\end{eqnarray}
\begin{eqnarray} 
\sum_{\mu=0}^\nu \left(\nu\atop \mu\right)[{\cal P}_\mu (x;l_1+l_2)
{\cal P}_{\nu-\mu}(x;0)-{\cal P}_\mu (x;l_1){\cal P}_{\nu-\mu}(x;l_2)]=0. 
\label{addtheor}
\end{eqnarray}
{\bf Proof.} The multiplicative structure of the generating function of 
Eq.~(\ref{3_8}) involving the parameter $l$ in the exponent implies the 
identity~(\ref{addid}). Substituting Eq.~(\ref{genfdef}) into this identity 
leads to Eq.~(\ref{addtheor}).  

We can also separate the $l$ dependence in the polynomials using 
Eq.~(\ref{sis}) in the generating function, Eq.~(\ref{3_6}). If 
$\sigma(x)=$~constant (as is the case for Hermite polynomials), the generating 
function only depends on the weight function, and the Taylor expansion of 
$w(x+y\sigma(x))$ for $|y\sigma(x)|<1$ is equivalent to the Rodrigues 
formula~(\ref{rods}).      

{\bf Corollary 1.}
\begin{eqnarray}\nonumber
{\cal P}(y,x;l)&=&\frac{w(x+y\sigma(x))}{w(x)}\left(1+y\sigma'(x)\right)^l\\
&=&\sum_{N=0}^{\infty}y^N\sum_{N-l\leq m\leq N}\left(l\atop N-m\right)
\sigma(x)^{m}\sigma'(x)^{N-m}\frac{w^{(m)}(x)}{m!w(x)}. 
\label{sep}
\end{eqnarray}
\begin{eqnarray}
{\cal P}_N(x;l)=\sum_{N-l\leq m\leq N}\left(l\atop N-m\right)
\frac{N!}{m!}{\cal P}_m(x;0)\sigma'(x)^{N-m}   
\label{sepex}
\end{eqnarray}
{\bf Proof.}  When $\sigma'(x)\neq 0,$ the Taylor expansion of the weight 
function in conjunction with a binomial expansion of the $l$th power of 
Eq.~(\ref{sis}) yields Eq.~(\ref{sep}). Using Eq.~(\ref{genfdef}) this 
translates into the polynomial expansion~(\ref{sepex}) that separates the 
complementary polynomials ${\cal P}_N(x;l)$ into the simpler polynomials 
${\cal P}_m(x;0)$ and the remaining $\sigma'(x)$ and $l$ dependence. 
Pearson's ODE~(\ref{wode}) guarantees the polynomial character of the 
${\cal P}_m(x;0)$ that are defined in Eq.~(\ref{rods}).  
 
Let us also mention the following symmetry relations.\\ 
{\bf Corollary 2.}
If $\sigma(-x)=(-1)^m\sigma(x), w(-x)=(-1)^n w(x)$ hold with integers $m,n$ 
then $P_l(-x)=(-1)^{l(m+1)}P_l(x)$ and  
\begin{eqnarray}
{\cal P}_{\nu}(-x;l)&=&{\cal P}_{\nu}(x;l),~m~{\rm odd},\label{symm1}\\ 
{\cal P}_{\nu}(-x;l)&=&(-1)^{\nu}{\cal P}_{\nu}(x;l),~m~{\rm even}.
\label{symm2}
\end{eqnarray}
{\bf Proof.}
The parity relation for $P_l(x)$ follows from substituting $-x$ in the 
Rodrigues formula~(\ref{rod}). The other polynomial parity relations follow 
from the identities  
\begin{eqnarray}
{\cal P}(y,-x;l)&=&{\cal P}(y,x;l),~m~{\rm odd}\label{symid1}\\
{\cal P}(-y,-x;l)&=&{\cal P}(y,x;l),~m~{\rm even}
\label{symid2}
\end{eqnarray}
which, in turn, result from substituting $-x$ into the first formula of 
Theorem 3. Expanding the generating functions according to their definition 
yields the relations~(\ref{symm1}), (\ref{symm2}).  

\section{Sturm--Liouville ODE}

{\bf Theorem 6.} The polynomials ${\cal P}_{\nu}(x;l)$ satisfy the 
Sturm-Liouville differential equation  
\begin{eqnarray}
\frac{d}{dx}\left(\sigma(x)^{l-\nu+1}w(x)\frac{d{\cal P}_\nu(x;l)}{dx}\right)=
-\lambda_\nu \sigma(x)^{l-\nu}w(x){\cal P}_\nu(x;l),
\label{sl}
\end{eqnarray}
which is equivalent to
\begin{eqnarray}
\sigma(x)\frac{d^2{\cal P}_\nu(x;l)}{dx^2}+[(l-\nu)\sigma'(x)+\tau(x)]
\frac{d{\cal P}_\nu(x;l)}{dx}=-\lambda_\nu {\cal P}_\nu(x;l),
\label{slode}
\end{eqnarray}
and the eigenvalues are given by 
\begin{eqnarray}
\lambda_{\nu}=-\nu [(l-\frac{\nu+1}{2})\sigma''+\tau'],~\nu=0, 1, \ldots  . 
\label{sleig}
\end{eqnarray}
{\bf Proof.} This is derived by natural induction again. The first step for 
$\nu=1$ is straightforward to verify. 

The step from $\nu$ to $\nu+1$ proceeds from the lhs of Eq.~(\ref{sl}) for 
$\nu+1,$ where we replace ${\cal P}_{\nu+1}$ by ${\cal P}_{\nu}$ using the 
recursive ODE~(\ref{rode}) so that, after some elementary manipulations, we 
end up with 
\begin{eqnarray}\nonumber
&&\frac{d}{dx}\left(\sigma(x)^{l-\nu}w(x)\frac{d{\cal P}_{\nu+1}(x;l)}{dx}
\right)=\sigma(x)^{l-\nu-1} w(x)\\\nonumber&\cdot&\bigg\{[(l-\nu-1)\sigma''
+\tau'][\sigma(x)\frac{d{\cal P}_\nu(x;l)}{dx}+[\tau(x)+(l-\nu-1)\sigma'(x)]
{\cal P}_\nu(x;l)]\\\nonumber&+&[(l-\nu-1\sigma'(x)+\tau(x))]
[\sigma(x)\frac{d^2{\cal P}_\nu(x;l)}{dx^2}+[(l-\nu)\sigma'(x)
+\tau(x)]\frac{d{\cal P}_\nu(x;l)}{dx}]\\\nonumber&+&\sigma(x)\frac{d}{dx}
[\sigma(x)\frac{d^2{\cal P}_\nu(x;l)}{dx^2}+[(l-\nu)\sigma'(x)+\tau(x)]
\frac{d{\cal P}_\nu(x;l)}{dx}]\bigg\}\\\nonumber&=& \sigma(x)^{l-\nu-1} w(x)
\{ [(l-\nu-1)\sigma''+\tau']{\cal P}_{\nu+1}(x;l)-\lambda_\nu 
{\cal P}_{\nu+1}(x;l)\}\\&=&-\lambda_{\nu+1}{\cal P}_{\nu+1}(x;l),   
\label{sld}
\end{eqnarray}
where we have used the recursive ODE~(\ref{rode}) and the ODE~(\ref{slode}) for
 the index $\nu$ repeatedly. Eq.~(\ref{rode}) introduces a third derivative of 
${\cal P}_\nu(x;l),$ a term which shows up as the next to last term on the rhs 
of the first equality sign in Eq.~(\ref{sld}). This completes the proof by 
induction and establishes the recursion 
\begin{eqnarray}
\lambda_{\nu+1}=\lambda_{\nu}-[(l-\nu-1)\sigma''+\tau']
\end{eqnarray}
for the eigenvalues, whose solution is Eq.~(\ref{sleig}). 

\section{Classical Polynomial Examples}

In the case of {\bf Hermite} 
polynomials~\cite{sz},\cite{aw},\cite{den},\cite{handbook},\cite{gr} $\sigma$ 
has no roots, so $\sigma(x)=$\\constant~$=1,$ without loss of generality, and 
$\sigma'=0;$ moreover, we may take $a_l=0, b_l=-2$ so $\tau(x)=-2x.$ Hence 
Pearson's ODE yields the weight function $w(x)=e^{-x^2}$ on $(-\infty,\infty)$ 
that is characteristic of Hermite polynomials. The Rodrigues 
formula~(\ref{rod}) then identifies the polynomials $P_l(x)=(-1)^l H_l(x)$ as 
Hermite's, while the Rodrigues formula~(\ref{rods}) for the complementary 
polynomials implies ${\cal P}_\nu(x;l)=P_\nu(x),$ so they are independent of 
the index $l$ and also Hermite polynomials. The recursive ODE~(\ref{rode}) 
becomes the well known differential recursion $H_{n+1}(x)=2xH_n(x)-H'_n(x).$ 
The Sturm-Liouville ODE becomes the usual ODE of the Hermite polynomials. The 
recursion~(\ref{recu1}) is the basic $H_{n+1}(x)=2x H_n(x)-2n H_{n-1}(x).$ 
Eq.~(\ref{brecu}) gives the differential recursion $H'_n(x)=2n H_{n-1}(x).$ 
The parity relation is also the well known one. The generating function is the 
standard one. Equation~(\ref{sep}) reproduces the usual expansion of Hermite 
polynomials in powers of the variable $x.$ 
 
For {\bf Laguerre} polynomials, $\sigma(x)$ has one real root, so $\sigma(x)=x$
 and $\tau(x)=1-x$ without loss of generality. Pearson's ODE gives the 
familiar weight function $w(x)=e^{-x}$ on $[0,\infty).$ Rodrigues 
formula~(\ref{rod}) identifies $P_l(x)=l!L_l(x).$ The Sturm-Liouville 
ODE~(\ref{slode}) 
\begin{eqnarray}
x\frac{d^2{\cal P}_\nu(x;l)}{dx^2}+(l+1-\nu-x)\frac{d{\cal P}_\nu(x;l)}{dx}
=-\lambda_\nu {\cal P}_\nu(x;l), \lambda_{\nu+1}=\lambda_\nu+1
\end{eqnarray}  
allows identifying ${\cal P}_\nu(x;l)=\nu!L_\nu^{l-\nu}(x)$ as an associated 
Laguerre polynomial. So, in the following we shift $l\to l+\nu,$ as a rule.   
The recursive ODE~(\ref{rode}) yields the differential recursion
\begin{eqnarray}
(\nu+1)L_{\nu+1}^{l-1}(x)=x\frac{d L_\nu^l(x)}{dx}+(l-x)L_\nu^l(x)
\end{eqnarray}   
which, in conjunction with  
\begin{eqnarray}
L_{\nu+1}^{l-1}(x)=L_{\nu+1}^{l}(x)-L_{\nu}^{l}(x),
\end{eqnarray}
leads to the standard three-term recursion
\begin{eqnarray}
(\nu+1)L_{\nu+1}^{l}(x)=(l+\nu+1-x)L_{\nu}^{l}(x)+x\frac{d L_{\nu}^{l}(x)}{dx}.
\end{eqnarray} 
The formula~(\ref{rods}) of {Theorem 2} is the usual Rodrigues formula for 
associated Laguerre polynomials, while the generalized Rodrigues 
formula~(\ref{b5}) 
\begin{eqnarray}
L_\nu^l(x)=\frac{\mu!}{\nu!}e^x x^{-l}\frac{d^{\nu-\mu}}{dx^{\nu-\mu}}\left(
x^{l+\nu-\mu} e^{-x} L_\mu^{l+\nu-\mu}(x)\right)
\end{eqnarray}  
is not part of the standard lore. 

The generating function~(\ref{3_8}) for this case becomes 
\begin{eqnarray}
L(y,x;l)=\sum_{\nu=0}^{\infty}y^\nu L_\nu^{l-\nu}(x)=e^{-xy}(1+y)^l
\end{eqnarray}
and is simpler than the usual one for associated Laguerre polynomials, which 
is the reason why our method is more elementary and faster than the standard 
approaches. The recursion~(\ref{recu1}) becomes 
\begin{eqnarray}
(\nu+1)L_{\nu+1}^{l-1}(x)=(l-x)L_{\nu}^l(x)-xL_{\nu-1}^{l+1}(x),
\end{eqnarray}
while the recursion~(\ref{recu2}) becomes 
\begin{eqnarray}
(\nu+1)L_{\nu+1}^l(x)=(l+\nu+1-x)L_{\nu}^l(x)-xL_{\nu-1}^{l+1}(x),
\end{eqnarray}
and Eq.~(\ref{recu3}) translates into 
\begin{eqnarray}
\frac{d L_{\nu-1}^{l+1}(x)}{dx}+\frac{d L_\nu^l(x)}{dx}=-L_{\nu-1}^{l+1}(x)
-L_{\nu-2}^{l+2}(x), 
\end{eqnarray}
a sum of the known recursion $\frac{d L_\nu^l(x)}{dx}=-L_{\nu-1}^{l+1}(x)$ 
which is the basic recursive ODE~(\ref{brecu}). Equation~(\ref{sep}) gives the 
standard expansion 
\begin{eqnarray}
L_N^l(x)=\sum_{n=0}^N\left(l+n\atop N-n\right)\frac{(-x)^n}{n!}. 
\end{eqnarray}

The simplest addition theorem originates from the elegant identity
\begin{eqnarray}
L(y,x_1;n_1)L(y,x_2;n_2)=L(y,x_1+x_2;n_1+n_2)
\end{eqnarray}
which translates into the polynomial addition theorem
\begin{eqnarray}
{\cal P}_\nu(x_1+x_2;n_1+n_2)=\sum_{k=0}^\nu\left(\nu\atop k\right)
{\cal P}_{\nu-k}(x_1;n_1){\cal P}_k(x_2;n_2)
\end{eqnarray}
and 
\begin{eqnarray}
L_\nu^{n_1+n_2}(x_1+x_2)=\sum_{k=0}^\nu L_k^{n_1-k}(x_1)
L_{\nu-k}^{n_2+k}(x_2)
\label{adf}
\end{eqnarray}
for associated Laguerre polynomials which is not listed in the standard 
ref.~\cite{gr} or elsewhere.

In the case of {\bf Jacobi} polynomials, $\sigma(x)$ has two real roots at 
$\pm1,$ without loss of generality; so 
\begin{eqnarray}
\sigma(x)=(1-x)(1+x),~\tau(x)=b-a-(2+a+b)x,
\end{eqnarray}
in a notation that will allow us to use the standard parameters.

Pearson's ODE~(\ref{wode}) leads to 
\begin{eqnarray}
w(x)=(1-x)^a (1+x)^b, 
\end{eqnarray} 
and Rodrigues formula~(\ref{rod}) and (\ref{rods}) identify the polynomials 
\begin{eqnarray}
P_l(x)=2^l(-1)^l l! P_l^{(a,b)}(x),~{\cal P}_\nu(x;l)=(-2)^\nu \nu!
P_\nu^{(a+l-\nu,b+l-\nu)}(x).
\end{eqnarray}
Thus, we shift $l\to l+\nu$  in translating our general results to Jacobi 
polynomials, as a rule. We may also set $l=0$ because this index merely shifts 
the parameters $a, b.$  

The recursive ODE~(\ref{rode}) translates into
\begin{eqnarray}\nonumber
-2(\nu+1)P_{\nu+1}^{(a-1,b-1)}(x)&=&[b-a-(a+b)x]P_{\nu}^{(a,b)}(x)\\&+&
(1-x^2)\frac{d P_{\nu}^{(a,b)}(x)}{dx}. 
\end{eqnarray}

The Sturm-Liouville ODE~(\ref{slode}) reproduces the usual ODE of Jacobi 
polynomials. The generating function, Eq.~(\ref{3_8}), 
\begin{eqnarray}\nonumber
{\cal P}(y,x;l)&=&\frac{[1-x-y(1-x^2)]^a[1+x+y(1-x^2)]^b}{(1-x)^a (1+x)^b}
\bigg\{\frac{1-(x+y(1-x^2)^2}{1-x^2}\bigg \}^l\\&=&[1-y(1+x)]^a[1+y(1-x)]^b
[1-2xy-y^2(1-x^2)]^l
\end{eqnarray}
is much simpler than the standard one~\cite{sz}, especially when we set $l=0,$ 
allowing for the transparent derivation of many recursion relations. For 
example, Eq.~(\ref{recu2}) becomes 
\begin{eqnarray}
-4(\nu+1)P_{\nu+1}^{(a-1,b-1)}(x)&=&2[b-a-x(a+b)]P_{\nu}^{(a,b)}(x)\\\nonumber 
&-&(1-x^2)[\nu+1+a+b]P_{\nu-1}^{(a+1,b+1)}(x),
\end{eqnarray}
Eq.~(\ref{recu2}) translates into
\begin{eqnarray}\nonumber
-4(\nu+1)P_{\nu+1}^{(a-1,b-1)}(x)&=&2[b-a-x(a+b+2\nu)]
P_{\nu}^{(a-1,b-1)}(x)\\&+&(1-x^2)[a+b+2\nu]P_{\nu}^{(a,b)}(x),
\end{eqnarray}
and Eq.~(\ref{recu3}) takes the form
\begin{eqnarray}\nonumber
&&(x^2-1)\frac{d P_{\nu-2}^{(a+2,b+2)}(x)}{dx}+4x\frac{d 
P_{\nu-1}^{(a+1,b+1)}(x)}{dx}+4\frac{d P_{\nu}^{(a,b)}(x)}{dx}\\
\nonumber&=& 2[a+b+2\nu]P_{\nu-1}^{(a+1,b+1)}(x)+
[b-a+x(a+b+2\nu)]\\&\cdot&P_{\nu-2}^{(a+2,b+2)}(x).
\end{eqnarray}
Equation~(\ref{sep}) gives
\begin{eqnarray}\nonumber
P_N^{(a, b)}(x)&=&(-2)^N N!(1-x)^{-a}(1+x)^{-b}\sum_{n=0}^N\left(N\atop n
\right)\\&\cdot&\frac{(-2x)^n(1-x^2)^{N-n}}{(N-n)!}\frac{d^n}{dx^n}
[(1-x)^a(1+x)^b]. 
\end{eqnarray}

A product formula for Jacobi polynomials is obtained from an addition
theorem in the variable $y$ for our generating function for $l=0$ (where we 
display the upper parameters now for clarity) 
\begin{eqnarray}\nonumber
&&{\cal P}^{(a,b)}(y_1,x;0){\cal P}^{(a,b)}(y_2,x;0)=
[(1+y_1(1-x)(1+y_2(1-x)]^b\\\nonumber&\cdot& 
[(1-y_1(1+x)(1-2(1+x)]^a\\\nonumber&=&
[1+(y_1+y_2)(1-x)]^b\{1+\frac{y_1y_2(1-x)^2}{1+(y_1+y_2)(1-x)}\}^b\\
\nonumber&\cdot&[1-(y_1+y_2)(1+x)]^a\{1+\frac{y_1y_2(1+x)^2}{1-(y_1+y_2)
(1+x)}\}^a\\\nonumber&=&{\cal P}^{(a,b)}(y_1+y_2,x;0)
\sum_{j,k=0}^{\infty}\left(a\atop k\right)\left(b\atop j\right)
\left(\frac{y_1y_2(1-x)^2}{1+(y_1+y_2)(1-x)}\right)^j\\\nonumber&\cdot&
\left(\frac{y_1y_2(1+x)^2}{1-(y_1+y_2)(1+x)}\right)^k=\sum_{j,k=0}^{\infty}
\left(a\atop k\right)\left(b\atop j\right)\\&\cdot&
{\cal P}^{(a-k,b-j)}(y_1+y_2,x;0)y_1^{j+k}y_2^{j+k}(1-x)^{2j}(1+x)^{2k}.  
\end{eqnarray}
Expanding into Jacobi polynomials according to Eq.~(\ref{genfdef}), 
comparing like powers of $y_1y_2,$ converting to Jacobi polynomials and 
shifting $a\to\\ a+\nu_1, b\to b+\nu_1$ yields the product formula
\begin{eqnarray}\nonumber
&&P_{\nu_1}^{(a,b)}(x)P_{\nu_2}^{(a+\nu_1-\nu_2,b+\nu_1-\nu_2)}(x)=
\sum_{0\leq\nu\leq (\nu_1+\nu_2)/2}2^{-2\nu}
\left(\nu_1+\nu_2-2\nu\atop \nu_1-\nu\right)\\\nonumber&\cdot&
\sum_{k=0}^\nu\left(a+\nu_1\atop k\right)\left(b+\nu_1\atop \nu-k
\right)(1+x)^{2k}(1-x)^{2(\nu-k)}\\&\cdot&
P^{(a+2\nu-\nu_2-k,b+\nu-\nu_2+k)}_{\nu_1+\nu_2-2\nu}(x). 
\label{prod}
\end{eqnarray}

\section{Conclusions}
We have used a natural way of working with the Rodrigues formula of a given 
set of orthogonal polynomials which leads to a set of closely related 
complementary polynomials that obey their own Rodrigues formulas, always have 
a generating function that can be summed in closed form leading to a 
transparent derivation of numerous recursion relations and addition theorems. 
These complementary polynomials satisfy a homogeneous second-order 
differential equation similar to that of the original polynomials.

Our method generates all the basics of the Hermite polynomials. It 
generates the {\bf associated} Laguerre polynomials and many of their known  
properties and new ones from the Laguerre polynomials in an elementary way.    
It also simplifies the derivations of various results for Jacobi polynomials. 

Our method is not restricted to the classical polynomials; when it is applied 
to the polynomials that are part of the wave functions of the Schr\"odinger 
equation with the Rosen-Morse and the Scarf potentials, it links these 
polynomials to the Romanovski polynomials which will be shown elsewhere. 



\begin{thebibliography}{0}

\bibitem{ni} A.\ F.\ Nikiforov and V.\ B.\ Uvarov, {\it Special Functions of 
Mathematical Physics,} Birkh\"auser Verlag, Basilea, 1988. 

\bibitem{is} M. E. H. Ismail, {\it Classical and Quantum Orthogonal 
Polynomials in One Variable,} Cambridge Univ. Press, Cambridge, 2005.

\bibitem{les} P. A. Lesky,  ``Endliche und unendliche Systeme von 
kontinuierlichen klassischen Orthogonalpoynomen,'' {\it Z.angew.Math.Mech.} 76 
(1996), 3, pp. 181-184.

\bibitem{sz} G. Szeg\"o, {\it Orthogonal Polynomials}, American Math. Soc., 
Vol. XXIII, Providence, RI, 1939.

\bibitem{aw}  G.\ B.\ Arfken and H.\ J.\ Weber, {\it Mathematical Methods for 
Physicists}, 6th ed., Elsevier-Academic Press, Amsterdam, 2005.

\bibitem{den} Phylippe Dennery and Andr\'e Krzywicki, {\it Mathematics for 
Physicists,} Dover, New York, 1996.

\bibitem{handbook}
M. Abramowitz and I. A. Stegun, {\it Handbook of Mathematical Functions with 
Formulas, Graphs and Mathematical Tables,} Dover, 2nd edition, New York, 1972.

\bibitem{gr} I. S. Gradshteyn and I. M. Ryzhik, {\it Table of Integrals, 
Series and Products}, ed. A. Jeffrey, Acad. Press, San Diego, 2000. 

\end{thebibliography}
\end{document}